\documentclass{amsart}
\usepackage{amsmath,amssymb}
\usepackage{color}
\usepackage{epsfig}
\usepackage{graphicx}
\usepackage[utf8]{inputenc}
\usepackage[T1]{fontenc}
 \usepackage[colorlinks=true]{hyperref}

\usepackage{dsfont}

\def\what{\widehat}
\def\RR{{\mathbb{R}}}
\def\ren{{\RR^N}}
\def\dyle{\displaystyle}

\def\e{\mathrm{e}}

\renewcommand{\d}{\, \mathrm{d}}

\newcommand{\pv}{\mathrm{P.V.}}
\newcommand{\lj}{\mathcal{L}}

\hypersetup{urlcolor=blue, citecolor=red}

\newtheorem{theorem}{Theorem}[section]
\newtheorem{corollary}[theorem]{Corollary}

\newtheorem{lemma}[theorem]{Lemma}
\newtheorem{proposition}[theorem]{Proposition}

\theoremstyle{definition}
\newtheorem{definition}[theorem]{Definition}
\newtheorem{remark}{Remark}

\title[Nonlocal heat equations] %Use the shortened version of the full title
      {Nonlocal heat equations: \\ decay estimates and Nash inequalities}

\author[Cristina Brändle and Arturo de Pablo]{}


 \email{cbrandle@math.uc3m.es}
 \email{arturop@math.uc3m.es}

\thanks{Work  partially supported by  Spanish project MTM2011-25287.}

\begin{document}

\newpage

\maketitle

% Enter the first author's name and address:
\centerline{\scshape C. Brändle and A. de Pablo}
\medskip
{\footnotesize
% please put the address of the first author
 \centerline{Departamento de Matem\'{a}ticas; Universidad Carlos III de Madrid }
   \centerline{Avda. de la Universidad 30}
   \centerline{28911 Leganés}
} % Do not forget to end the {\footnotesize by the sign }

\begin{abstract}
We obtain  $L^q$--$L^p$ decay estimates, $1\le q<p<\infty$ for solutions of nonlocal heat equations of the form $\partial_tu+\lj u=0$. Here $\lj$ is an integral operator given by a symmetric nonnegative kernel of L\'evy type. We obtain these estimates in terms only of the behaviour of the kernel at infinity, without any information of its behaviour at the origin. This includes bounded and unbounded transition probability densities. An equivalence between the decay and a restricted Nash inequality is shown.  We  also prove that $\lim_{t\to \infty}\|u(t)\|_\infty=0$.
Finally we deal with  nonlinear nonlocal equations of porous medium type $\partial_tu+\lj\Phi(u)=0$.
\end{abstract}

\bigskip

%%%%%%%%%%%%%%%%%%%%%%%%%%%%%%%%%%%%%%%%%%%%%%%%
\section{Introduction}\label{sect-introduction}
\setcounter{equation}{0}
The purpose of this work is to study the decay for large time  of  solutions to nonlocal evolution equations  of the form
\begin{equation}\label{NLHE}\partial_tu+\lj u=0,\quad x\in\ren,\;t>0.
\end{equation}
Here the diffusion is driven by an integral operator $\lj$ in the spatial variable   defined  by
\begin{equation}\label{NLHOp}
  \lj f(x)=\pv\int_{\ren}(f(x)-f(y))J(x-y)\d y,
\end{equation}
with $J$ a Lévy kernel, i.e., satisfying $
\int_{\ren}J(z)(|z|^2\wedge1)\d z<\infty.
$
When studying the long-time behaviour of solutions for problems like \eqref{NLHE}--\eqref{NLHOp}, it is common in the literature to impose a global behaviour on $J$, typically integrability on the whole space, or boundedness, or a power-type behaviour, which may be different at the origin and at infinity, see for instance \cite{BassLevin2002,CarlenKusuokaStroock1987,ChasseigneChavesRossi2006,IgnatRossi2009}. As we will see this is it not always necessary: in this paper we obtain estimates on the decay of solutions for large times in terms only of the behaviour of the kernel at  infinity,  without imposing any condition at the other end. To this purpose we establish some functional inequalities that might have independent interest.

Problem~\eqref{NLHE}--\eqref{NLHOp} has been widely studied in  many different contexts, see for instance \cite{BarlowBassChenKassmann2009,BassLevin2002,KassmannSchwab2013,Mimica2012,Sztonyk2011}. Typically we find in the literature two lines of research, often disconnected, namely $J$ integrable or $J$ hypersingular of fractional Laplacian type. These two kinds of kernels give rise to very different properties and characteristics both in the nature of the problem and the solutions.

In particular, regarding the operator, it is of differential type or not depending on the singularity of the kernel at the origin. For regular kernels or even weakly singular kernels, $J\in L^1(\ren)$, and to fix ideas say $\int_{\ren}J=1$, the operator in~\eqref{NLHOp} can be written in the form
\begin{equation*}
  \label{convolution}
  \lj f(x)=f(x)\int_{\ren}J(x-y)\d y-\int_{\ren}f(y)J(x-y)\d y=f(x)-J*f(x),
\end{equation*}
where $*$ denotes the usual convolution in $\mathbb{R}^N$. This means that $\lj f$ is as regular as $f$, and $\lj$ is a zero order operator.
On the other hand, if $J(z)=|z|^{-N-\alpha}$, with $0<\alpha<2$, then the condition $f\in C^{\alpha+\varepsilon}(\ren)$, with $0<\alpha+\varepsilon<1$, implies $\lj f\in C^\varepsilon(\ren)$, see~\cite{Silvestre2007} where more general cases are also considered. The operator $\lj$ is in that case (a multiple of) the well known \emph{fractional Laplacian} $(-\Delta)^{\alpha/2}$, an operator of order $\alpha$, like $\alpha$ derivatives. For general L\'evy kernels we get that the operator $\lj$ is well defined for smooth bounded functions and it is of order strictly less than 2, see for instance Lemma~\ref{lem:c-alpha}.

\medskip

\noindent{\sc About the solution -- }
The operator $\lj$ given by formula~\eqref{NLHOp} satisfies for $f$ in the Schwartz class,
\begin{equation}
  \label{lj-symbol-form}
 \widehat{(\lj f)}(\xi)=m(\xi)\widehat f(\xi),
\end{equation}
for some multiplier $m(\xi)$, cf. formula~\eqref{symbol-L}, it is therefore a \emph{pseudodifferential operator}, see~\cite{Stein1993}.
The advantage of formula \eqref{lj-symbol-form} is that equation~\eqref{NLHE} becomes $\partial_t\what u+m\what u=0$. If we are given an initial value $u(x,0)=u_0(x)$, then we get the formal expression of the solution,
\begin{equation}
 \label{eq:u.sol.prob}
u(x,t)=u_0*\mu_t(x),\qquad \what{\mu_t}(\xi)=\e^{-m(\xi)t}.
\end{equation}
(We also write $u(t)=u_0*\mu_t$ by abuse of notation). This solution will be as regular as  the most regular of the two terms in the convolution. In fact, the  \emph{fundamental solution} $\mu_t$ may be a regular function or may be a singular measure depending on the behaviour of the multiplier.  We will devote Subsection~\ref{sect:multiplier}
to study the properties of  $m$ in terms of the Lévy kernel.

\medskip

\noindent{\sc Regularizing effect --}
 The first issue is then to study the regularity properties of $\mu_t$ for $t$ small. We characterize this regularity using the behaviour of the multiplier at the origin.
Regarding again the two different lines of research mentioned above, if $J$ is integrable there is no regularizing effect, as $\what\mu_t$ is not integrable and then $\mu_t$ is discontinuous; actually $\mu_t$ is a singular measure, see for instance \cite{ChasseigneChavesRossi2006}. On the very other hand, for the fractional Laplacian type operators, the solution is immediately bounded and $C^\infty$. We want to determine which is the threshold for the kernel $J$ in order to get regularizing effect or not. In the brief Subsection~\ref{subsect:regular} we prove that this borderline is precisely the function $J(z)=|z|^{-N}$ for $|z|$ small, and it lies into the no regularizing effect side, though it is not integrable, see Theorem~\ref{thm:main1}.

\medskip

\noindent{\sc Decay estimates --}
We devote the main part of the paper to study the decay of the solution for large times imposing  the least as possible conditions on the kernel $J$, in particular regarding only the behaviour  of the kernel at infinity. Decay estimates in $L^p$ for  singular kernels at the origin of the type of the fractional Laplacian have been obtained in \cite{BassLevin2002,CarlenKusuokaStroock1987,Komatsu1988}. On the other hand, using in a definite way that the kernel is bounded and for $N\ge2$, some decay estimates are proved in~\cite{ChasseigneFelmerRossiTopp2013}, see also~\cite{IgnatRossi2009}.

In our case, we obtain decay estimates without imposing any condition on the kernel at the origin: the kernel could be bounded, integrable or even singular of fractional Laplacian type. In other words, the associated transition probability densities could  be bounded or unbounded.
This includes for instance the borderline case $J(z)\sim|z|^{-N}$ for $|z|$ small, which is not covered by the previous studies. This generality makes that the specific arguments used before do not work.

Assuming that the kernel posses what is known as a {\it heavy tail}, that is, it decays slower than a power $|z|^{-N-\alpha}$, for large $|z|$ and some $0<\alpha<2$ (infinite second order momentum), we prove for some $t$ large and $1\leq q<p<\infty$, the estimate
\begin{equation}\label{fract-decay}
\|u(t)\|_p\le ct^{-\varrho(q,p)}\|u_0\|_q, \qquad \varrho(q,p)=\frac{N}{\alpha}\left(\frac1q-\frac1p\right),
\end{equation}
see Theorem~\ref{th:decay-asymp-sharp}. We remark that this result holds for every $N\ge1$. If no condition is imposed on $J$ at infinity we always get the same estimate, namely~\eqref{fract-decay} but with $\alpha=2$.

\medskip

\noindent{\sc Nash inequalities --} Decay estimates for heat like equations, or more generally, for submarkovian strongly continuous symmetric semigroups, are usually shown to be equivalent  to some Sobolev or Nash inequalities, see for instance~\cite{CarlenKusuokaStroock1987, Nash1958, Varopoulos1985}. Here we establish some relation between the large time decay for this kind of semigroups with a restricted Nash inequality for the associated Dirichlet form, see Section~\ref{sect:nash}. We also prove that  this Nash inequality indeed holds in our context of equation~\eqref{NLHE}, see Corollary~\ref{cor-inversa-H}.

\medskip

\noindent{\sc Nonlinear problem --}
We are also interested in nonlinear nonlocal heat equations of the form
\begin{equation*}\label{NLNLHE}
  \left\{
  \begin{array}{rcl@{\qquad}l}
    \partial_tu+\lj\Phi(u)&=&0,&\quad x\in\ren,\;t>0, \\
    u(x,0)&=&u_0(x),&\quad x\in\ren,
  \end{array}
  \right.
\end{equation*}
for some constitutive function $\Phi$, which behaves (in a weak sense) like a power $\sigma\ge1$ for small values. In Subsection~\ref{subsect:nonlinear} we show that the procedure used to deal with the linear problem works also for this nonlinear case. We obtain again  $L^q$--$L^p$ decay estimates similar as before, but under the  limitation that $\max\{1,\sigma-1\}<q<p<\infty$, see Theorem~\ref{th:asymp-nonlinear}.

\medskip

\noindent{\sc Organization of the paper --} In the preliminary Section~\ref{sect:preliminaries} we settle the problem and the functional framework in which we will work throughout the paper. Section~\ref{sect:nash} is self-contained and relates the long-time behaviour of the solutions of an abstract problem with some functional inequalities of Nash type. The long Section~\ref{sect:behaviour} contains the main results of the paper, namely the decay estimates for large times of the solutions to equation~\eqref{NLHE}. This is done in Subsection~\ref{subsect:decay}. We also characterize the regularizing effect, Subsection~\ref{subsect:regular}. For both studies we use the properties of the associated multiplier, described in~Subsection~\ref{sect:multiplier}. In Subsection~\ref{subsect:nonlinear} we apply the previous method to study the decay for the solutions of a nonlinear nonlocal Porous Medium type equation.

\section{Preliminaries. Problem setting}
\label{sect:preliminaries}

We introduce in this section different notions of solution and establish some basic properties for  the Cauchy problem associated to equation~\eqref{NLHE},
\begin{equation}\label{NLHE-problem}
\left\{
  \begin{array}{rcl@{\qquad}l}
    \partial_tu+\lj u&=&0,&\quad x \in\ren,\;t>0, \\
   u(x,0)&=&u_0(x),&\quad x\in\ren.
  \end{array}
  \right.
\end{equation}
 We will be assuming throughout the paper (without mentioning it anymore) that $J$ is nonnegative, radially symmetric, positive near the origin and  satisfies the L\'evy condition
\begin{equation*}\label{levy}
\int_{\ren}J(z)(|z|^2\wedge1)\d z<\infty.
\end{equation*}
The \lq\lq classical\rq\rq\   examples mentioned in the Introduction, $J$ integrable and $J=|z|^{-N-\alpha}$, are in particular Lévy  kernels (provided $0<\alpha<2$) and they somehow delimit the range of admissible kernels.

  In order to give a sense to~\eqref{eq:u.sol.prob}, we can consider classical solutions, for which we require that the initial data are bounded $C^{1,1}$ functions. We  also construct weak solutions for $L^1(\mathbb{R}^N)\cap L^2(\mathbb{R}^N)$ initial values. Weaker conditions on $u_0$ can also be considered depending on $\lj$. For operators $\lj$ given by bounded kernels, an existence and uniqueness theory is established in~\cite{ChasseigneChavesRossi2006} for $u_0\in L^1(\mathbb{R}^N)$ and such that $\what u_0\in L^1(\mathbb{R}^N)$.
We present some results concerning the basic properties of the solutions to problem~\eqref{NLHE-problem}. The proofs are straightforward and we omit them.

Using the symmetry of $J$ we can write~\eqref{NLHOp} in the  equivalent formulation
\begin{equation}\label{NLHOp2}
\lj f(x)=\frac12\int_{\ren}(2f(x)-f(x+z)-f(x-z))J(z)\d z,
\end{equation}
which allows to prove that the operator $\lj$ is well defined for smooth enough functions.
\begin{lemma}\label{lem:c-alpha}
If $f\in C^{1,1}(\ren)\cap L^\infty(\ren)$ then $\lj f\in C(\ren)\cap L^\infty(\ren)$.
\end{lemma}

If moreover $f$ is compactly supported then, for $|x|$ large it is easy to prove that $
|\lj f(x)|\leq c\sup_{|z|>|x|/2} J(z)
$.
Thus, we consider the weighted  space ${\mathcal F}_\lj:=L^1(\mathbb{R}^N,\rho\d x)$, where $\rho(x)=\min\{1,\sup_{|z|>|x|/2} J(z)\}$.
Thus the operator $\lj$ can be defined in the sense of distributions for functions $f\in {\mathcal F}_\lj$. For instance, if $J$ is nonincreasing this allows to define $\lj$ for locally integrable functions increasing with some small growth at infinity depending on $J$.

Observe that the operator $\lj$ can be defined also in the sense of distributions for functions in a larger space. For instance, if $J$ is nonincreasing this allows to define $\lj$ for locally integrable functions increasing at infinity less than $|x|^N$.

\medskip

\noindent {\sc Weak solutions --} As it is common in the literature, we may define a \emph{weak solution} of problem~\eqref{NLHE-problem} by formally multiplying the equation by a suitable test function and then integrating by parts.
In our situation this is done by considering  the quadratic Dirichlet form (nonlocal interaction energy)
\begin{equation}
  \label{quadratic-form}
 \mathcal E(f,f)=\frac12\int_{\ren}\int_{\ren}(f(x)-f(y))^2J(x-y)\d x\d y,
\end{equation}
and the associated bilinear form
\begin{equation}
  \label{bilinear-form}
 \mathcal E(f,g)=\frac12\int_{\ren}\int_{\ren}(f(x)-f(y))(g(x)-g(y))J(x-y)\d x\d y.
\end{equation}
Then we define
\begin{equation}
  \label{lj-quadratic-form}
 \langle \lj f,g\rangle=\mathcal E(f,g),
\end{equation}
for functions in the space where the quadratic form is finite. Actually we consider the space
\begin{equation*}
  \label{def-HL}
  \mathcal{H}_{\lj}(\RR^N)=\{f\in L^2(\RR^N)\,:\,\mathcal E(f,f)<\infty\},
\end{equation*}
with associated norm
\begin{equation*}
  \label{norm-HL}
  \|f\|_{\mathcal{H}_{\lj}}^2=\|f\|_{2}^2+ \mathcal E(f,f).
\end{equation*}
We also use the equivalent expression of the Dirichlet form in terms of the multiplier,
\begin{equation}
  \label{quadratic-form-mult}
 \mathcal E(f,f)=\int_{\ren}m(\xi)|\widehat f(\xi)|^2\d\xi,
\end{equation}
see \eqref{lj-symbol-form}.
It is clear that the space $\mathcal{H}_{\lj}(\ren)$ contains  $H^1(\ren)$. Roughly speaking, it is of Sobolev type if $\lj$ is of positive order. In fact, the more singular is $J$ at the origin, the smaller is $\mathcal{H}_{\lj}(\ren)$. On the other side, if $J\in L^1(\mathbb{R}^N)$ then $\mathcal{H}_{\lj}(\ren)\equiv L^2(\mathbb{R}^N)$ and  $\|f\|_{\mathcal{H}_{\lj}}^2\leq c\|f\|_{2}^2\|J\|_{1}$.

\begin{definition}
    Let $u_0\in L^1(\mathbb{R}^N)\cap L^2(\mathbb{R}^N)$.
    A weak solution of problem~\eqref{NLHE-problem} is a function $u\in  L^1(\mathbb{R}^N\times(0,T))\cap L^2((0,T); \mathcal{H}_{\lj}(\mathbb{R}^N))$, for every $T>0$, such that
\begin{equation}
  \label{eq:def.weak}
   \int_0^\infty \int_{\mathbb{R}^N} u\,\partial_t\varphi\d x\d t+ \int_{\mathbb{R}^N} \varphi(x,0)u_0(x)\d x=\int_0^\infty\mathcal{E}(u(t),\varphi(t))\d t
   \end{equation}
    for every $\varphi\in C^1(\mathbb{R}^N\times(0,\infty))$, vanishing for $|x|$ and $t$ large.
\end{definition}

\begin{remark}\label{rem-test}  As we have said, if $\varphi\in C^1_0(\mathbb{R}^N)$ then $\mathcal E(\varphi,\varphi)<\infty$ and then the last integral makes sense. In fact more general test functions can be used in the definition, for instance $\varphi\in L^2((0,\infty);\mathcal{H}_{\lj}(\mathbb{R}^N))$, such that $\partial_t\varphi\in L^2(\mathbb{R}^N\times(0,T))$, $\varphi(\cdot,t)=0$ for $t\ge T$.
\end{remark}
%
%We recall the function defined in the introduction

%which is a solution to the problem in Fourier variables
%
%
Problem~\eqref{NLHE-problem} can be written, whenever the Fourier Transform is well defined, using~\eqref{lj-symbol-form} as
\begin{equation*}\label{eq:prob-fourier}
\left\{
  \begin{array}{rcl@{\qquad}l}
    \partial_t\widehat u+m\widehat u&=&0, \\
    \widehat u(\xi, 0)&=&\widehat u_0(\xi).
  \end{array}
  \right.
\end{equation*}
This ODE problem has the explicit solution $\widehat u(\xi,t)= \widehat u_0(\xi)\e^{-m(\xi)t}$.
The fact that this problem is equivalent to our original problem in weak form is the content of the next theorem and in particular leads to an explicit solution of~\eqref{NLHE-problem}.

\begin{theorem}\label{th:u.sol}
  For every initial datum $u_0\in L^1(\mathbb{R}^N)\cap L^2(\mathbb{R}^N)$, the function \begin{equation}
 \label{eq:u.sol.prob2}
u(t)=u_0*\mu_t,\qquad \what{\mu_t}(\xi)=\e^{-m(\xi)t},
\end{equation}
  is the  unique weak solution to problem~\eqref{NLHE-problem}.
\end{theorem}

The next properties imply that the weak solution is in fact a \emph{strong solution}, and equation~\eqref{NLHE-problem} holds almost everywhere. Moreover, in view of Remark~\ref{rem-test}, we can use the solution $u$ itself as a test function. This reveals to be of great utility in proving the decay of the solution for large times, Section~\ref{sect:behaviour}.

\begin{lemma}
  \label{lemma:properties}
  Let $u$ given by~\eqref{eq:u.sol.prob2} be the weak solution to~\eqref{NLHE-problem}.
  \begin{itemize}
    \item If $u_0\in L^p(\ren)$, $1\le p\le\infty$, then $u(t)\in L^p(\ren)$ for every $t>0$,  and $\|u(t)\|_p\le\|u_0\|_p$. Even more, if $u_0\ge0$ then $\int_{\ren}u(x,t)\d x=\int_{\ren}u_0(x)\d x$.
    \item $u(t)\in \mathcal{H}_{\lj}(\mathbb{R}^N)$ for every $t>0$, and $\|u(t)\|_{\mathcal{H}_{\lj}}\le (1+1/(2\e t))^{1/2}\|u_0\|_2$.
    \item $\partial_tu(t)\in L^2(\mathbb{R}^N)$ for every $t>0$.
  \end{itemize}
\end{lemma}

  \medskip

\noindent{\sc Very weak solutions --}
We consider now an even weaker concept of solution, using the weighted space $\mathcal{F}_\lj$.
\begin{definition}
    Let $u_0\in\mathcal{F}_\lj$.
    A very weak solution of problem~\eqref{NLHE-problem} is a function $u\in  L^1((0,T);\mathcal{F}_\lj)$, for every $T>0$, such that
\begin{equation}
  \label{eq:def.very.weak}
   \int_0^\infty \int_{\mathbb{R}^N} u\,\partial_t\varphi\d x\d t+ \int_{\mathbb{R}^N} \varphi(x,0)u_0(x)\d x=\int_0^\infty \langle u,\mathcal{L}\varphi\rangle\d t
   \end{equation}
    for every $\varphi\in C^2(\mathbb{R}^N\times(0,\infty))$, vanishing for $|x|$ and $t$ large.
\end{definition}

By approximation with data in $L^1(\mathbb{R}^N)\cap L^2(\mathbb{R}^N)$ we get existence of very weak solutions.    Uniqueness  is not clear. For a similar concept of solution in the case of the fractional Laplacian operator, in~{\rm\cite{BarriosPeralSoriaValdinoci2013}} uniqueness is proved by using strongly the boundedness of $\mu_t$ for $t$ positive.

\begin{lemma} Let  $u_0\in \mathcal{F}_\lj$ be an initial data for~\eqref{NLHE-problem}. Then $u$  given in~\eqref{eq:u.sol.prob2}
  is a very weak solution to problem~\eqref{NLHE-problem}.
\end{lemma}

  \medskip

{\noindent{\sc Classical solutions --}
Lemma~\ref{lem:c-alpha} says that $\lj u$ is continuous in space provided $u$ is $C^{1,1}$ and bounded in space for all $t>0$. On the other hand, regularity in time follows from the fact that the $L^2$-function $\partial_t\widehat{u}(t)=-m\e^{-mt}\widehat{u}_0$ is continuous in time. Hence, the weak solution $u$ to problem~\eqref{NLHE-problem} is in that case a \emph{classical solution}, and the equation~\eqref{NLHE-problem} holds pointwisely.

Imposing additional conditions on $J$ we get regularity in space, {\it a posteriori}, under very few conditions on the initial data.
For instance, if $J(z)\le c|z|^{-N}$ for $|z|$ small then, exactly as in Lemma~\ref{lem:c-alpha}, $u_0\in C^\varepsilon(\mathbb{R}^N)\cap L^\infty(\mathbb{R}^N)$ implies $\lj u(t)\in C(\mathbb{R}^N)$ and the solution is again classical. When $J\in L^1(\mathbb{R^N})$, then it is enough to have $u_0$ continuous and bounded in order to have a classical solution. Finally, as we will see in  Section~\ref{sect:behaviour}, when $J(z)\gg |z|^{-N}$ for $|z|$ small, then the function $\widehat\mu_t$ decays faster than any power for each fixed positive time, and so $\mu_t$ is $C^\infty$ and bounded for every $t>0$. This gives a classical solution for every initial value $u_0\in L^1(\mathbb{R^N})$.

\section{Nash inequalities}\label{sect:nash}

In this section we  relate the decay of the solutions to problem~\eqref{NLHE-problem} with some functional inequalities of Nash type. The application of the Nash inequalities to get the precise decay estimates will be  the subject of Section~\ref{subsect:decay}.

As a starting point we   recall a result  due to Varopoulos, see \cite{Varopoulos1985}, that characterizes the decay of  heat semigroups $\mathcal{T}_t=e^{-At}$, for positive self-adjoint operators $A$ in $L^2$, in terms of functional inequalities. Actually, if $\mathcal{T}_t$  satisfies the \emph{ultracontractivity property}
\begin{equation}
  \label{decay-general-heat}
  \|\mathcal{T}_t(v)\|_\infty\le ct^{-d}\|v\|_1
\end{equation}
for some $d>1$ and any $t>0$, then a Sobolev type inequality
\begin{equation}
 \label{Sob-type}
 \mathcal{Q}(f,f)\ge c\|f\|_{\frac{2d}{d-1}}^2
 \end{equation}
 holds for the Dirichlet form $\mathcal{Q}(f,f)=\langle Af,f\rangle$ and  any function $f$ in the domain $\mathfrak{D}$ of $\mathcal{Q}$. The converse is also true.
In particular, in the case of the Laplacian, and even in the fractional Laplacian case $A=(-\Delta)^{\alpha/2}$, $0<\alpha\le2$, estimate~\eqref{decay-general-heat} is easy to obtain. Here $d=N/\alpha$, from the homogeneity of the semigroup. When $N>\alpha$ this   gives  the well-known \emph{Hardy-Littlewood-Sobolev inequality}
\begin{equation}
 \label{HLS}
 \|(-\Delta)^{\alpha/4}f\|_{2}\ge c\|f\|_{\frac{2N}{N-\alpha}}.
 \end{equation}

The proof of~\eqref{decay-general-heat} from~\eqref{Sob-type} uses that~\eqref{Sob-type} implies another interesting inequality~\eqref{nash-cks}, called \emph{Nash inequality}, see~\cite{Nash1958}. Observe that~\eqref{nash-cks} makes sense for every $d>0$, not only for $d>1$. The equivalence is shown in~\cite{CarlenKusuokaStroock1987}.

\begin{theorem}[\cite{Nash1958, CarlenKusuokaStroock1987}]
  \label{th:cks}
  Let $\mathcal Q$ and $\mathcal{T}_t$ be the Dirichlet form and the heat semigroup, resp.,  associated to a positive self-adjoint operator  in $L^2$, and let $d>0$.  The following conditions are equivalent:
  \begin{equation}
    \label{nash-cks}
    \mathcal Q(f,f)\ge C_1\|f\|_{2}^{2+2/d}\qquad\forall\;f\in \mathfrak{D}\cap L^1(\mathbb{R}^N),\quad \|f\|_1=1,
  \end{equation}
for some (for every) $1<p\le\infty$,
\begin{equation}
    \label{hypercontractivity}
    \|\mathcal{T}_t(v)\|_p\le C_2t^{-d(1-1/p)}\|v\|_1\qquad\forall\;t>0,\quad v\in L^1(\mathbb{R}^N).
  \end{equation}
\end{theorem}

Later \cite{Coulhon1996} extend Theorem~\ref{th:cks} above  writing \eqref{nash-cks} in the form
\begin{equation}
    \label{nash-davies}
    \mathcal Q(f,f)\ge C_1\|f\|_{2}^{2}B(\|f\|_{2}^{2})\qquad\forall\;f\in \mathfrak{D}\cap L^1(\mathbb{R}^N),\quad \|f\|_1=1,
  \end{equation}
where $B(s)$ can be any function in a certain class generalizing the function $B(s)=s^{1/d}$ of~\eqref{nash-cks}. They obtain an estimate of the semigroup
\begin{equation}
    \label{hypercon-2}
    \|\mathcal{T}_t(v)\|_\infty\le C_2\Lambda(t)\|v\|_1\qquad\forall\;t>0,\quad v\in L^1(\mathbb{R}^N),
  \end{equation}
for some function $\Lambda$ related to $B$. Actually, the admissible functions $B$ always satisfy the estimate $B(s)\gg \log s$ for $s\to\infty$.

As we will see in Section~\ref{subsect:regular}, a decay like \eqref{hypercon-2} for the solution of~\eqref{NLHE-problem} cannot be true for any function $\Lambda$ and all times, if no assumption of singularity at the origin of the kernel $J$ defining $\lj$ is made.
Hence it is natural to ask if there is a  Nash inequality of the type \eqref{nash-davies} so that  a decay of the form \eqref{hypercon-2} can be obtained only for large times.  We show that a result of this kind is true with the function $B(s)=\min\{1,s^{1/d}\}$. We remark that this function $B$ is bounded. Unfortunately we have been not able to reach the value $p=\infty$.

\begin{theorem}
  \label{th:nash->decay}
  Let $\mathcal Q$ and $\mathcal{T}_t$ be as before. Assume there exist constants $1\le r<2$, $d>0$ and $C_1>0$ such that every function $f\in\mathfrak{D}\cap L^1(\mathbb{R}^N)\cap L^2(\mathbb{R}^N)$ with $\|f\|_r=1$ satisfies
  \begin{equation}
    \label{our-nash}
    \mathcal Q(f,f)\ge C_1\|f\|_{2}^{2}\min\{1,\|f\|_{2}^{2/d}\}.
  \end{equation}
  Then, for any $1\le q<p<\infty$ and every $v\in L^1(\mathbb{R}^N)\cap L^p(\mathbb{R}^N)$,  there exist a constant $C_2>0$ and a time $t_0>0$ such that
  \begin{equation}
    \label{our-hypercon}
    \|\mathcal{T}_t(v)\|_p\le C_2t^{-\varrho(q,p)}\|v\|_q\qquad\forall\;t>t_0,\qquad \varrho(q,p)=\frac{dr}{2-r}\left(\frac1q-\frac1p\right).
  \end{equation}
\end{theorem}

\begin{remark}  If $\|f\|_2$ is small, i.e., $\|f\|_2\leq 1$, and $r=1$, inequality~\eqref{our-nash} reduces to the Nash inequality~\eqref{nash-cks}; otherwise it is a Poincaré inequality, $\mathcal{Q}(f,f)\ge c\|f\|_2^2$.
 \end{remark}

The proof uses the following inequality due to Stroock \cite{Stroock1984} and Varopoulos \cite{Varopoulos1985}.
\begin{proposition}[\cite{Stroock1984, Varopoulos1985}]
  \label{pro:SV-original} For any nonnegative function $f\in \mathfrak{D}$ such that $f^a,\,f^b\in \mathfrak{D}$, with $a,\,b\ge0$, and $a+b=2$, it holds
  \begin{equation}\label{SV-original}
\mathcal{Q}(f^a,f^b) \ge ab
\mathcal{Q}(f,f).
\end{equation}
\end{proposition}

{\noindent{\it Proof of Theorem} \ref{th:nash->decay}. }
Let $q\ge1$ be given and take $p=2q/r$. Put $z(t)=\mathcal{T}_t(v)$.
We have
\begin{equation}\label{semingroup-weak}
\int_{\mathbb{R}^N}\partial_tz\zeta=-\mathcal{Q}(z,\zeta)
\end{equation}
for every test function $\zeta$. Putting $\zeta=|z|^{p-2}z$, and using \eqref{SV-original}, we get
$$
\frac{\d}{\d t}\int_{\mathbb{R}^N}|z|^p\le -c\mathcal{Q}(|z|^{p/2},|z|^{p/2}).
$$
We now apply the hypothesis~\eqref{our-nash} to the function $|z|^{p/2}$ to get
$$
\mathcal{Q}(|z|^{p/2},|z|^{p/2}) \ge c \|z\|_p^{p}\min \left\{1,\, \left(\frac{\|z\|_p} {\|z\|_q}\right)^{p/d}\right\}.
$$
This holds for every $t>0$. Therefore, denoting $\psi(t)=\|z(t)\|_p^p$, since the $L^q$-norms are nonincreasing in time, see~\cite{Fukushima1980}, that is, $\|z(t)\|_q\le \|z(0)\|_q=\|v\|_q $, we obtain the differential inequality
$$
\psi'(t)+
c\psi(t)\min\left\{1,\,\psi(t)^{1/d}\|v\|_q^{-p/d}\right\}\le0.
$$
If  $\psi(t)\geq \|v\|_q^p$ for every $t>0$, the previous inequality gives $\psi'(t)\leq -c\psi(t)$, which implies $\lim_{t\to\infty}\psi(t)=0$ and this is a contradiction.
Hence,  there exists $0\le t_1<\infty$ such that
$\psi(t)\le \|v\|_q^p$ for $t=t_1$, and hence for any $t>t_1$. Therefore $\psi$ satisfies
\begin{equation}
\label{eq:inequality.decay}
\left\{
\begin{array}{ll}
\psi'(t)\le-c\|v\|_q^{-p/d}\psi(t)^{1+1/d},&\quad t>t_1,\\ [3mm]
\psi(t_1)\le \|v\|_q^p.
\end{array}
\right.
\end{equation}
The time $t_1$ depends on the initial data. Indeed,
$
  t_1=cp\log_+\left( \frac{\|u_0\|_p}{\|u_0\|_q}\right).
$
We solve the differential inequality to get the estimate
$$\psi(t)\le \|v\|_q^p\left(1+\frac{c}{d}(t-t_1)\right)^{-d}, \qquad t>t_1,
$$
and finally, taking $t_0=2t_1$,
\begin{equation}
  \label{eq:starting.point.iteration}
\|z(t)\|_p\le c\|v\|_qt^{-d/p}, \qquad t>t_0.
\end{equation}
Observe that $d/p=dr/2q=\varrho(q,2q/r)$.

The general case follows easily by iteration and interpolation. Indeed,
if we define $p_k=(2/r)^kq$ we get
$$
\|z(t)\|_{p_k}\leq c_k t^{-d/p_k}\|z(t_{k-1}\|_{p_{k-1}}, $$
for $t-t_{k-1}$ large. This gives
$$
\|z(t)\|_{p_k}\leq c t^{-\sum_{j=1}^k d/p_j}\|v\|_q=c t^{-\varrho(q,p_k)}\|v\|_q,
$$
for $t$ large enough.
Finally, for an arbitrary $q<p<\infty$ we use interpolation with some $p_k>p$ and the fact that $L^q$-norm decays.
\qed

We now prove the following converse result.
\begin{theorem}
  \label{th:main-converse}
 Let $1\le q<2$ and $p= q/(q-1)$ ($p=\infty$ if $q=1$), and $\tau>0$, $\nu>0$. For any function $v\in\mathfrak{D}\cap L^1(\mathbb{R}^N)\cap L^p(\mathbb{R}^N)$ satisfying
 \begin{equation}
    \label{our-decay}
    \|\mathcal{T}_t(v)\|_p\le C_2t^{-\nu}\|v\|_q\qquad\forall\;t>\tau,
  \end{equation}
  there holds
    \begin{equation}
    \label{our-nash2}
    \mathcal Q(v,v)\ge C\|v\|_{2}^{2}\min\left\{\frac{1}{\tau},\left(\frac{\|v\|_{2}}{\|v\|_{q}}\right)^{2/\nu}\right\}.
  \end{equation}
\end{theorem}

\proof Put $z(t)=\mathcal{T}_t(v)$. We first have
$$\int z(t)v\leq \|z(t)\|_{p}\|v\|_q\leq ct^{-\nu}\|v\|_q^2.$$
On the other hand, following for instance~\cite{CarlenKusuokaStroock1987}, choosing $\zeta=v$ in \eqref{semingroup-weak}, we get
$$
\int z(t)v= \|v\|_2^2-\int_0^t \mathcal{Q}(z(s),z)\geq \|v\|_2^2-t\mathcal{Q}(v,v),
$$
since clearly, by the definition of $\mathcal{Q}$ we have $\mathcal{Q}(z(s),z(0))\leq \mathcal{Q}(z(0),z(0))=\mathcal{Q}(v,v)$ for every $s>0$.
Summing up
$$
\|v\|_2^2\leq t\mathcal{Q}(v,v)+c t^{-\nu}\|v\|_q^2\qquad\mbox{for every } t>\tau.
$$
We minimize in $t$ the function $h(t)=At+Bt^{-\nu}$. If  $t_1=(\nu B/A)^{1/(\nu+1)}>\tau$, then
$$
\|v\|_2^2\leq c \|v\|_q^{2/(\nu+1)} \mathcal{Q}(v,v)^{\nu/(\nu+1)}.
$$
On the contrary, the estimate obtained is
$$\|v\|_2^2\leq \tau\mathcal{Q}(v,v)+ \tau^{-\nu}\mathcal{Q}(v,v)\tau^{\nu+1} .$$
Hence
$$
\mathcal{Q}(v,v)\geq C\|v\|_2^2\min \left\{\frac{1}{\tau}, \|v\|_2^{2/\nu}\|v\|_q^{-2/\nu}\right\}.
$$
\qed

We remark that by Hölder's inequality the result can be extended to any $p>q/(q-1)$ if $q>1$. The exponent $\nu$ in \eqref{our-nash2} must be replaced by $\mu=\frac{(2-q)p\nu}{p-q}$.

\section{Behaviour of the solutions}\label{sect:behaviour}

\subsection{Properties of the multiplier}\label{sect:multiplier}

Formulas~\eqref{lj-symbol-form} and \eqref{eq:u.sol.prob} show  that the operator $\mathcal{L}$ and the solution to problem~\eqref{NLHE-problem} are characterized only in terms of the multiplier (or symbol in the probabilistic context) $m$. This multiplier  has the well known Lévy-Khintchine expression, see for instance~\cite{Bertoin1996},
\begin{equation}
  \label{symbol-L}
m(\xi)=\int_{\mathbb{R}^N} \big (1-\cos(z\cdot \xi)\big)J(z)\d z.
\end{equation}
Since  the kernel $J$ is of L\'evy type we have that $m$ is well defined and moreover we obtain the first estimate
\begin{equation}
  \label{eq:lower.bound.m}
C_1\min\{1,|\xi|^2\}\le m(\xi)\le C_2\max\{1,|\xi|^2\}.
\end{equation}

Of course, a better knowledge of the behaviour of $J$ would imply sharper bounds for $m$. For instance, if $J$ is integrable then $m(\xi)\le\|J\|_{1}$.
On the other hand, if $J$ has a fractional Laplacian type singularity at the origin, say $J(z)\sim|z|^{-N-\alpha}$ for some $0<\alpha<2$, then $m(\xi)\sim|\xi|^\alpha$ for $|\xi|$ large.
In general, if we define
\begin{equation}
  \label{def.ell}
  \ell(r)=|z|^NJ(z), \qquad r=|z|,
\end{equation}
then  the behaviour of $\ell$ at the origin has a direct translation on the behaviour $m$ for large values of $|\xi|$. In fact, these relation has already been studied in~\cite{KassmannMimica2013} in a different context related to Hölder estimates of an associated elliptic problem. We consider for a given function $\ell$ the following functions,
$$
  \psi_1^\ell(r)=\int_{r}^1\frac{\ell(s)}s\d s,\qquad
  \psi_2^\ell(r)=r^{-2}\int_0^r s\ell(s)\d s.
  $$
As we have said before, we must consider only non-integrable kernels and thus for $\ell $ as in~\eqref{def.ell}, $\lim_{r\to 0} \psi_1^\ell(r)=\infty$. On the other hand, the integral in $ \psi_2^\ell$ is well defined since $J$ is a Lévy kernel.

Now, the upper estime for the multiplier is straightforward using these functions.
\begin{lemma}
\label{lemma:m.for.xi.large.up} The multiplier $m(\xi)$ satisfies, for $|\xi|>1$,
  \begin{equation*}
 m(\xi)\le c_1\psi_1^\ell(1/|\xi|)+c_2\psi_2^\ell(1/|\xi|),
  \end{equation*}
  \end{lemma}
A lower bound is obtained in~\cite{KassmannMimica2013} with the extra hypothesis that
\begin{equation}
  \label{eq.ell.beta}\ell(r)\geq \beta(r), \qquad \text{for $r\in(0,1)$},
\end{equation}
 where $\beta$ is a positive function that
varies regularly at zero with index $\rho\in(-2,0]$;  that is
$$
\lim_{r\to 0^+}\frac{\beta(\lambda r)}{\beta(r)}=\lambda^\rho, \qquad \text{for every $\lambda>0$},
$$
see~\cite{BinghamGoldieTeugels1989}.

\begin{lemma}
\label{lemma:m.for.xi.large.low} Let $\ell$ satisfy~\eqref{def.ell} and~\eqref{eq.ell.beta} where $\beta$ is a regular varying function at zero. Then the multiplier $m(\xi)$ satisfies, for $|\xi|>1$,
  \begin{equation*}
   m(\xi)\ge c\psi_1^\beta(1/|\xi|).
  \end{equation*}
  \end{lemma}
The key point of the proof consists in estimating the integral in \eqref{symbol-L} only over the intervals where $1-\cos(z\cdot \xi)\ge 1/2$ and then using the properties of regular varying functions, see~\cite{KassmannMimica2013} for the details. We also notice that for if $\ell$ is a regular varying function, then it always holds $\psi^\ell_2\le c\psi^\ell_1$ near the origin, cf.~\cite{BinghamGoldieTeugels1989}. Therefore we have in that situation $m(\xi)\sim\psi^\ell_1(1/|\xi|)$ for $|\xi|$ large.

 We are now concerned with estimates for $m$ when $\xi$ is small. If we have a precise power-type lower bound of the decay at infinity of $J$ then the bound~\eqref{eq:lower.bound.m}  can be improved.
\begin{lemma}
  \label{lem:multip}
Assume  $J(z)\ge c|z|^{-N-\alpha}$ for $z$ large and some $\alpha>0$. Then $m(\xi)$ satisfies
\begin{equation}
    \label{multip}
    m(\xi)\ge c|\xi|^\gamma\qquad \text{for $|\xi|\le1$ and $\gamma=\min\{\alpha,2\}$}.
  \end{equation}
\end{lemma}

\proof
As we have said, if $\alpha\ge 2$ then~\eqref{multip} holds with $\gamma=2$. Let then
$\alpha\in(0,2)$. For $|\xi|$ small enough,
$$
\begin{aligned}
m(\xi)&\displaystyle=\int_{\mathbb{R}^N}(1-\cos(z\cdot \xi))J(z)\d z\geq c\int_{c_1|\xi|^{-1}<|z|<c_2|\xi|^{-1}}|\xi|^2|z|^2|z|^{-N-\alpha}\d z  \\
&\displaystyle= c|\xi|^2\int_{c_1|\xi|^{-1}}^{c_2|\xi|^{-1}}r^{1-\alpha}\d r\ge  c|\xi|^{\alpha}.
\end{aligned}$$
\qed

\subsection{Short-time behaviour}\label{subsect:regular}

As it is shown in Theorem \ref{th:u.sol}, the solution to problem~\eqref{NLHE-problem} is given by the convolution expression~\eqref{eq:u.sol.prob2}. Thus, if  $\mu_t$ is bounded, then problem~\eqref{NLHE-problem} has an $L^1$--$L^\infty$ regularizing effect in the form  $\|u(t)\|_\infty\le C(t)\|u_0\|_1$ for any $t>0$. Moreover, if $\mu_t$ is regular so is the weak solution even if the initial value is not. The fact that $\mu_t$ is bounded or not can be characterized by the integrability of $\widehat\mu_t$, and this last depends on the behaviour of the multiplier $m(\xi)$ at infinity. This behaviour also determines the regularity of the solution. In turns, all these properties depend only on the singularity of the kernel $J$ at the origin, through formula~\eqref{symbol-L}.

The model case for the regularizing effect is the fractional Laplacian kernel, $J(z)=|z|^{-N-\alpha}$, with $\alpha\in(0,2)$, for which the following estimate holds
\begin{equation}
\|u(t)\|_\infty\le ct^{-N/\alpha}\|u_0\|_1,
\label{smoothing-FHE}\end{equation}
for every $t>0$. In fact, since the multiplier is $m(\xi)=|\xi|^\alpha$ then $\what \mu_t=\e^{-|\xi|^\alpha t}$ and $\|\mu_t\|_\infty=\|\widehat\mu_t\|_1=ct^{-N/\alpha}$. The same estimate \eqref{smoothing-FHE} is valid, but only for bounded times, if the kernel satisfies $J(z)\le c|z|^{-N-\alpha}$ for $|z|<\eta$ and some $\eta>0$ small, since then, from Lemma~\ref{lemma:m.for.xi.large.low} we have
\begin{equation}\label{dosdecays}
\|\mu_t\|_\infty\le \int_{|\xi|<1}\e^{-c|\xi|^2t}\d\xi+\int_{|\xi|>1}\e^{-c|\xi|^\alpha t}\d\xi+\le c_1t^{-N/2}+c_2t^{-N/\alpha}.
\end{equation}
The decay of the solution for large times is at least that of the (local) heat equation, and further estimates depend on the behaviour of $J$ at infinity, see for instance~\cite{CarlenKusuokaStroock1987} and next subsection.

On the other hand, if $J$ is integrable in all of $\mathbb{R}^N$, and say $\|J\|_1=1$, then $m(\xi)=1-\what J(\xi)$. Hence, the Riemann-Lebesgue Lemma implies  $\lim_{|\xi|\to\infty}m(\xi)=1$ and  therefore $\what \mu_t\not\in L^p(\ren)$ for any $t>0$ and any $p<\infty$. Thus $\mu_t\not\in L^\infty(\mathbb{R}^N)$, and it is even discontinuous. This means that problem~\eqref{NLHE-problem} has no regularizing effect. In fact, an explicit expression for the convolution (singular in this case) measure  $\mu_t$ is given in \cite{ChasseigneChavesRossi2006}, where $\mu_t = \e^{-t} \delta_0 + \omega_t$, with $\delta_0$ the Dirac delta measure at the origin and $\omega_t$ as smooth as $J$. Therefore the solutions are at most as regular as the initial data are.

Next we characterize the borderline to have $L^1$--$L^\infty$ short time regularizing effect; that is, when $\what \mu_t\in L^1(\ren)$ for $t$ small. We show that the threshold for  having short time regularizing effect is not given by the integrability or not of $J$ at the origin, but in some sense by the function $J(z)=|z|^{-N}$.

\begin{theorem}\label{thm:main1}
Let $\ell$ be given in~\eqref{def.ell}.
\begin{itemize}
  \item If $\ell(r)\le c$ then there is no regularizing effect for problem~\eqref{NLHE-problem} for any small time.
  \item If $\ell$ varies regularly near zero, then there is  regularizing effect for small times for problem~\eqref{NLHE-problem} if and only if  $\lim_{r\to0}\ell(r)=\infty$.
\end{itemize}
\end{theorem}

\noindent{\it Proof. }
The condition $\ell(r)\le c$ implies  $\psi^\ell_1(r)\le c\log 1/r$ and $\psi^\ell_2(r)\le c$. Thus
$$
\int_{\ren}\what \mu_t(\xi)\d\xi\ge c\int_{|\xi|>1}\e^{-ct\log|\xi|}\d\xi=c\int_{|\xi|>1} \frac1{|\xi|^{ct}}\d\xi=\infty
$$
for every $t$ small.

On the contrary if $\lim_{r\to0}\ell(r)=\infty$, using Lemma~\ref{lemma:m.for.xi.large.low},
$$
\infty= \lim_{r\to0}\ell(r)=-c\lim_{r\to0}r(\psi^\ell_1)'(r)=c\lim_{r\to0}\frac{\psi^\ell_1(r)}{\log(1/r)}\leq \lim_{|\xi|\to\infty}\frac{m(\xi)}{\log|\xi|}.
$$
This gives that for every $t>0$ there exists $M$ large such that
$$
\begin{aligned}
\int_{\ren}\what \mu_t(\xi)\d\xi&\le c_1+\int_{|\xi|>M}\e^{-tm(\xi)}\d\xi\le c_1+\int_{|\xi|>M} \e^{-(N+1)\log|\xi|}\d\xi \\ &=c_1+\int_{|\xi|>M} \frac1{|\xi|^{N+1}}\d\xi<\infty.
\end{aligned}
$$
\qed

The condition $\ell(r)$ regularly varying cannot be avoided in the above Theorem. Indeed, we can construct singular oscillating functions $\ell(r)$ such that no regularizing effect occur. This shows that in fact the regularizing effect  depends on the \lq\lq measure'' of the singular part. Let
$$
\ell(r)=
\begin{cases}
  b_k, & \text{ if } a_{2k+1}<r<a_{2k},\\
  1, & \text{ if } a_{2k+2}<r<a_{2k+1},
\end{cases}
$$
with   $a_{2k}=2^{-k}$, $a_{2k+1}=2^{-k}(1-b_k^{-1})$, and any sequence $b_k$ of positive numbers diverging to infinity. Then there is not regularizing effect no matter what the sequence $b_k$ is.
We just observe that for $r\in(2^{-k},2^{-(k-1)})$,
$$
\begin{array}{rl}
\psi_1^\ell(r)&\displaystyle\leq c\sum_{j=0}^{k+1}\left(b_j\log\Big(\frac{a_{2j}}{a_{2j+1}}\Big)+\log\Big(\frac{a_{2j+1}}{a_{2j+2}}\Big)\right) \\ [4mm]&\displaystyle\leq c\sum_{j=0}^{k+1}\left(b_j\log\Big(\frac{1}{1-b_j^{-1}}\Big)+2\log2\right)\leq  ck,
\end{array}$$
as well as
$$
\begin{array}{rl}
\psi_2^\ell(r)&\displaystyle\leq c2^{2k}\sum_{j=k-1}^\infty\left(b_j\Big(a_{2j}^2-a_{2j+1}^2\Big)+\Big(a_{2j+1}^2-a_{2j+2}^2\Big)\right) \\ [4mm]&\displaystyle\leq c2^{2k}\sum_{j=k-1}^\infty 2^{-2j}\leq  c.
\end{array}$$
Thus
$$
m(\xi)\le ck\le c\log|\xi| \quad\text{for } |\xi|\sim 2^{k},
$$
and as before $\widehat\mu_t\not\in L^1(\mathbb{R}^N)$. Observe also that taking for instance $b_k=2^{\alpha k}$, we have no regularizing effect for a very singular kernel, $\limsup_{z\to 0} |z|^{N+\alpha}J(z)\geq c>0$.

\

The same type of estimates as the one used in Theorem~\ref{thm:main1} allow to determine the regularity of the fundamental solution $\mu_t$, and thus that of the solution. We summarize the regularizing effect in the following.
\begin{corollary} Let $u$ be a weak solution to problem~\eqref{NLHE-problem} with $u_0\in L^1(\mathbb{R}^N)$, $u_0\not\in L^\infty(\mathbb{R}^N)$.
\begin{itemize}
  \item If $J(z)\ge c|z|^{-N}$ for small $|z|$ then there exists $t_0>0$ such that $u(t)\in L^\infty(\mathbb{R}^N)$ for any $t>t_0$. Moreover, for any positive integer $k$ there exists $t_1>t_0$ such that $u(t)\in C^k(\mathbb{R}^N)$ for $t>t_1$.
    \item If $J(z)\gg |z|^{-N}$ for small $|z|$ then $u(t)\in L^\infty(\mathbb{R}^N)\cap C^\infty(\mathbb{R}^N)$ for any $t>0$.
    \item If $J(z)\le c|z|^{-N}$ for small $|z|$ then there exists $t_2>0$ such that $u(t)\not\in L^\infty(\mathbb{R}^N)$ for any $0<t<t_2$. If in addition $J\in L^1(\mathbb{R}^N)$ then $t_2=\infty$.
  \end{itemize}
\end{corollary}

\begin{remark}
  There exist kernels $J\not\in L^1(\mathbb{R})^N$ for which there is no regularizing effect for any positive time. For instance consider for $0<p\leq 1$, $J(z)=|z|^{-N}(\log(|z|^{-1}))^{-p}$.

  We also want to stress that in the elliptic framework,~\cite{KassmannMimica2013} show  a regularity property for $\mathcal{L}$-harmonic functions  for every nonintegrable Lévy kernel; i.e there is a different threshold to have regularity for the solution in the elliptic and parabolic problems.
\end{remark}

\subsection{Long-time behaviour}\label{subsect:decay}

When there is no regularizing effect, or in general, when~\eqref{decay-general-heat} does not hold for all $t>0$, but only for $t$ large, inequalities~\eqref{Sob-type} and~\eqref{nash-cks} cannot be true.
However, we expect that the restricted Nash inequality~\eqref{our-nash} must hold, reflecting the possible nonsingular behaviour of the kernel at the origin.
Since our aim now is proving the decay of the solutions for $t$ large, we concentrate in proving that~\eqref{our-nash} is true in the context of the space $\mathcal{H}_{\lj}(\mathbb{R}^N)$, with the Dirichlet form $\mathcal{E}(v,v)$. The symbol $m$ coming from a Lévy kernel $J$ satisfies~\eqref{multip}.

\begin{lemma}
\label{lemma:sobolev}
  Let $m$ satisfy \eqref{multip} for some $0<\gamma\le2$. Then,   for every  $1<r<s\le2$,
    \begin{equation}\label{mas-que-sobolev}
\|z\|_s^2\le c_1\|z\|_r^{2\theta_1}\mathcal{E}(z,z)^{1-\theta_1}+c_2\|z\|_r^{2\theta_2}\mathcal{E}(z,z)^{1-\theta_2},
  \end{equation}
where $c_i=c_i(r,s,\gamma,N)$ and
\begin{equation}\label{mas-que-sobolev-theta}
\theta_1=\dfrac{r[N(2-s)+\gamma s]}{s[N(2-r)+\gamma r]},\qquad
 \theta_2=\dfrac{r(2-s)}{s(2-r)}.
  \end{equation}\end{lemma}

The case $s=2$, $\gamma=2$ and $N\ge3$ is obtained in \cite{IgnatRossi2009}. They also are able to reach $r=1$ in~\eqref{mas-que-sobolev} provided $J$ is integrable. We mention that more general kernels $J=J(x,y)$ not only $J=J(x-y)$ are considered in that paper.

\proof Assume first $N>\gamma$. As in \cite{IgnatRossi2009} and \cite{Nash1958} we write $z=v+w$ with $\widehat v=\chi_{Q}\widehat z$, $Q=\{|\xi_k|\leq 1,\, k=1,\cdots N\}$ and $w=z-v$. It is clear that
$$
\mathcal{E}(z,z)=\mathcal{E}(v,v)+\mathcal{E}(w,w),
$$
and both terms on the right-hand side are positive. Using that $m$ verifies~\eqref{multip}, and by the Hardy-Littlewood-Sobolev inequality~\eqref{HLS}, we find the following bound for the first term,
$$
  \mathcal{E}(v,v)\geq c\int_{\mathbb{R}^N}|\xi|^\gamma|\widehat v(\xi)|^2\d\xi= c\|(-\Delta)^{\gamma/4}v\|_{2}^2\ge c\|v\|_{\frac{2N}{N-\gamma}}^2.
$$
Since $\|v\|_{r}\leq c\|z\|_{r}$ for every $r>1$, see~\cite{Stein1993} (the result is false for $r=1$), using  interpolation we get
\begin{equation}
  \label{eq:1}
   \|v\|_{s}
\leq \|v\|^{\theta_1}_{r}\|v\|^{1-\theta_1}_{{\frac{2N}{N-\gamma}}}\le c \|z\|^{\theta_1}_{r}\|\mathcal{E}(z,z)^{\frac{1-\theta_1}2},
\end{equation}
where $\frac{1}{s}=\frac{\theta_1}{r}+\frac{(1-\theta_1)(N-\gamma)}{2N}$.
As to the second term, the calculus is easier. We have
$$
  \mathcal{E}(w,w)\geq c\int_{\mathbb{R}^N}|\widehat w(\xi)|^2\d\xi= c\|w\|_{2}^2.
$$
So again by interpolation, and since $\|w\|_{r}\le\|z\|_{r}+\|v\|_{r}\leq c\|z\|_{r}$,
\begin{equation}
  \label{eq:2}
\|w\|_{s}\leq\|w\|_{r}^{\theta_2}\|w\|_{2}^{1-\theta_2}\leq c\|z\|_{r}^{\theta_2}\mathcal{E}(z,z)^{\frac{1-\theta_2}2},
\end{equation}
where $\frac1s=\frac{\theta_2}r+\frac{1-\theta_2}2$.
Adding inequalities~\eqref{eq:1} and~\eqref{eq:2} we get~\eqref{mas-que-sobolev}.

If $ N\in (\gamma/2, \gamma]$, i.e $N=\gamma=2$ or $N=1\leq \gamma<2$,
we cannot use  Hardy-Littlewood-Sobolev inequality. However, we may overcome this problem by defining
the new Dirichlet form
$$
\mathcal{E}_{1/2}(f,f)=\int_{\mathbb{R}^N} |\widehat f(\xi)|^2 m^{1/2}(\xi)\d \xi,
$$
and replacing $\gamma$ by $\gamma/2$ in the above argument. Observe that  the critical Sobolev exponent reads now $4N/(2N-\gamma)$. We then have
$$
\|v\|_s\le c_1\|v\|_r^{\theta'_1}\mathcal{E}_{1/2}(v,v)^{\frac{1-\theta'_1}2},
$$
where $\theta_1$, given in~\eqref{mas-que-sobolev-theta}, is replaced  here by  the new value $\theta'_1=\frac{r[2N(2-s)+\gamma s]}{s[2N(2-r)+\gamma r]}$, and
$$
\|w\|_s\le c_2\|z\|_r^{\theta_2}\mathcal{E}(z,z)^{\frac{1-\theta_2}2},
$$
where
$\theta_2$ remains the same. In order to write the expression for $v$ in terms of $\mathcal{E}(z,z)$ we use Hölder inequality,
$$
\mathcal{E}_{1/2}(v,v)=\int_{\mathbb{R}^N} m^{1/2}(\xi)|\widehat v(\xi)|^2 \d\xi
\leq \|\widehat v\|_2\left(\int_{\mathbb{R}^N}m(\xi) |\widehat v(\xi)|^2 \d\xi\right)^{1/2}=\|v\|_2\mathcal{E}(v,v)^{1/2}$$
and interpolation
$$
\|v\|_2\le\|v\|_r^a \|v\|_{\frac{2N}{N-\gamma/2}}^{1-a}\le \|v\|_r^a\mathcal{E}_{1/2}(v,v)^{\frac{1-a}2},
$$
with $a=\frac{r\gamma}{2N(2-r)+\gamma r}$,
to get
$$
\mathcal{E}_{1/2}(v,v)\le c\Big(\|v\|_r^{2a}\mathcal{E}(v,v)\Big)^{1/(1+a)}.
$$
Therefore we conclude
 $$
 \begin{aligned}
 \|z\|_s&\le c_1\|z\|_r^{\theta'+\frac{a(1-\theta_1')}{1+a}}\mathcal{E}(z,z)^{\frac{1-\theta_1'}{2(1+a)}}
 + c_2\|z\|_r^{\theta_2}\mathcal{E}(z,z)^{\frac{1-\theta_2}{2}} \\
 &=c_1\|z\|_r^{\theta_1}\mathcal{E}(z,z)^{\frac{1-\theta_1}{2}}
 + c_2\|z\|_r^{\theta_2}\mathcal{E}(z,z)^{\frac{1-\theta_2}{2}}.
 \end{aligned}
$$

Finally,  if $N=1$, $\gamma=2$, we define $\mathcal{E}_{1/4}$ and replace $\gamma$ by $\gamma/4$ in the same way as before. We omit the details.
~\qed

As a consequence of this lemma, we prove now that inequality~\eqref{our-nash} holds in our setting. In fact we state a more general version with the new exponent $1<s\le2$ (there $s=2$ and hence $\theta_2=0$). The interest of this generality comes from its applicability to a nonlinear problem of porous medium type, see Section~\ref{subsect:nonlinear}. Moreover, the result that we prove here holds for every $1<r<2$, not just for some $r$. The case $r=1$ is left open.

\begin{corollary}\label{cor-inversa-H} Let $m$ satisfy \eqref{multip} for some $0<\gamma\le2$ and fix  the parameters $r$ and $s$ such that $1<r<s\le2$. It holds
    \begin{equation*}\label{inversa-H}
\mathcal{E}(z,z) \ge c \min \left\{  \left(\|z\|_s\|z\|_r^{-\theta_1}\right)^{2/(1-\theta_1)},
\,\left(\|z\|_s\|z\|_r^{-\theta_2}\right)^{2/(1-\theta_2)} \right\},
 \end{equation*}
where $\theta_1$ and $\theta_2$ are given in~\eqref{mas-que-sobolev-theta}.
\end{corollary}

\proof
We define $X:=\mathcal{E}(z,z)^{1/2}$,  $Y:=\|z\|_s$ and $K:=\|z\|_r$. Inequality~\eqref{mas-que-sobolev} can then be written as
$$
Y\leq c_1 K^{\theta_1} X^{1-\theta_1}+ c_2K^{\theta_2} X^{1-\theta_2}\leq\left\{\begin{array}
  {l@{\qquad}l}
  \displaystyle c_3K^{\theta_1} X^{1-\theta_1}, &X\leq K,\\[8pt]
  \displaystyle c_3K^{\theta_2} X^{1-\theta_2},& X\geq K,
\end{array} \right.
$$
$c_3=c_1+c_2$. Hence, $X\geq  (Y/(c_3K^{\theta_1}))^{1/(1-\theta_1)}$ if $X\leq K$ and $X\ge (Y/(c_3K^{\theta_2}))^{1/(1-\theta_2)}$ otherwise.
That is
$
X\ge c\min\big\{K^{-\theta_1/(1-\theta_1)}Y^{1/(1-\theta_1)},K^{-\theta_2/(1-\theta_2)}Y^{1/(1-\theta_1)}\big\}.
$~\qed

We are now able to prove one of the main results of the paper, namely the $L^q$--$L^p$ estimates for large times for the solutions to problem~\eqref{NLHE-problem}.
As we have said, it depends on the behaviour of the kernel $J$ at infinity.
If we look at the estimates of the multiplier~\eqref{eq:lower.bound.m} and~\eqref{multip}, it is clear that imposing no condition at infinity on $J$ will give the same result as imposing the decay condition $J(z)\ge c|z|^{-N-\alpha}$ for large $|z|$ with $\alpha=2$. We therefore only consider the case $0<\alpha\le2$.
\begin{theorem}
  \label{th:decay-asymp-sharp}
  Let $1<p<\infty$ and $\bar p=\max\{p,2\}$. Assume $J(z)\ge c|z|^{-N-\alpha}$ for large $|z|$ with $0<\alpha\le2$. Given $u_0\in L^1(\mathbb{R}^N)\cap L^{\bar p}(\mathbb{R}^N)$, let $u$ be the associated  solution to problem~\eqref{NLHE-problem}. Then there exists $t_0>0$ large such that, for every $t>t_0$ and $1\le q<p$,
  \begin{equation}
    \label{decay-p-q}
    \|u(t)\|_p\le ct^{-\frac{N}\alpha\left(\frac1q-\frac1p\right)}\|u_0\|_q.
  \end{equation}
\end{theorem}

\proof
  It follows directly from Lemma~\ref{lem:multip}, Corollary~\ref{cor-inversa-H} with $s=2$ and any $r\in(1,2)$, and Theorem~\ref{th:nash->decay} with $d=\frac{N(2-r)}{r\alpha}$.~\qed

If we  try to find an ultracontractivity type estimate for the $L^\infty$-norm for large times, we face several difficulties that we address.
First of all, the time $t_0$ for which the estimate~\eqref{decay-p-q} holds, depends on $p$ and it tends to infinity as $p$ grows. Hence, passing to the limit in the proof of Theorem~\ref{th:decay-asymp-sharp} or Theorem~\ref{th:nash->decay} is not possible.
On the other hand, the classical duality argument fails. Indeed, we cannot consider the heat semigroup associated to the operator $\lj$ at a fixed time $t$, $S_t:L^1\to L^p$, since again the time $t$ for which it is a bounded operator depends on the $L^p$-norm of the solution itself.

Of course, if we impose  additional conditions to the kernel $J$ or the initial condition $u_0$, then some decay of the solution can be proved. For instance, if $J(z)\ge c|z|^{-N-\beta}$ for $|z|$ small, $0<\beta<2$, then $
\|u(t)\|_\infty\le C(t^{-N/\beta}+t^{-N/\alpha})\|u_0\|_1,
$ see~\eqref{dosdecays} and
Corollary {\rm 2.12} in~\cite{CarlenKusuokaStroock1987}.

On the other hand, if the initial datum $u_0$ in problem~\eqref{NLHE-problem} satisfies the very restrictive condition $\what u_0\in L^1(\mathbb{R}^N)$ (in particular it requires  $u_0$  continuous), then a decay estimate for the solution is easy to obtain. Actually,
$$
\begin{aligned}
  \|u(t)\|_\infty&\le \int_{|\xi|<1}\e^{-c_1|\xi|^\alpha t}|\what u_0(\xi)|\d\xi+\int_{|\xi|>1}\e^{-c_2t}|\what u_0(\xi)|\d\xi \\
  &\le c_3\|u_0\|_1t^{-N/\alpha}+c_4\|\what u_0\|_1\e^{-c_2 t}.
\end{aligned}
$$

In view of this, our purpose here is  to obtain at least some behaviour for large times assuming only that the initial data is $u_0\in L^1(\mathbb{R}^N)\cap L^\infty(\mathbb{R}^N)$, and of course with no condition on $J$ at the origin.

\begin{theorem}
  \label{th:decay-asymp-sharp-infty}
  Let  $u_0\in L^1(\mathbb{R}^N)\cap L^\infty(\mathbb{R}^N)$ and assume that the hypotheses of Theorem~{\rm\ref{th:decay-asymp-sharp}} hold. Then
  \begin{equation}
    \label{decay-infty}
    \lim_{t\to\infty }\|u(t)\|_\infty=0.
  \end{equation}
\end{theorem}

\proof
Let $S_t$ be the heat semigroup associated to $\lj$. Theorem~\ref{th:decay-asymp-sharp} with for instance $q=1$ and $p=2$, implies that
$$
\lim_{t\to\infty }t^{\varrho}\|S_t(u_0)\|_2\leq C,\qquad \varrho=\frac{N}{2\alpha}.
$$
Thus  there is a subsequence of times $t_k\to\infty$ such that
$\lim_{k\to\infty}t_k^{\varrho}S_{t_k}(u_0)$ converges a.e. to a bounded function in $L^2(\mathbb{R}^N)$.
Hence, since $S_t$ is selfadjoint,
$$
\begin{aligned}
\|\lim_{k\to\infty}t_k^{\varrho}S_{t_k}(u_0)\|_\infty
  &=\sup_{\|v\|_1\leq 1}\left|\int_{\mathbb{R}^N}
\lim_{k\to\infty}t_k^{\varrho}S_{t_k}(u_0)v\right|\\
&\leq \sup_{\|v\|_1\leq 1} \lim_{k\to\infty}t_k^{\varrho}\left|\int_{\mathbb{R}^N} u_0S_{t_k}(v)\right|
\\
&\leq \sup_{\|v\|_1\leq 1}\lim_{k\to\infty} t_k^{\varrho}\|S_{t_k}(v)\|_2 \|u_0\|_2\leq c\|u_0\|_2.
\end{aligned}
$$
This means that
$$
\lim_{k\to\infty}t_k^{\varrho}\|u(t_k)\|_\infty\leq c\|u_0\|_2,$$
and \eqref{decay-infty} follows since the $L^\infty$-norm is nonincreasing.
\qed

\subsection{Decay estimates for a nonlinear equation}\label{subsect:nonlinear}

We want to apply the  method described in Subsection~\ref{subsect:decay} to deal with the nonlocal nonlinear filtration equation
\begin{equation}\label{NLNLHE2}
  \left\{
  \begin{array}{rcl@{\qquad}l}
    \partial_tu+\lj\Phi(u)&=&0,&\quad x\in\ren,\;t>0, \\
    u(x,0)&=&u_0(x),&\quad x\in\ren.
  \end{array}
  \right.
\end{equation}
As to the constitutive nonlinear function we assume it is of Porous Medium type in weak sense, that is, $\Phi\in C^1(\mathbb{R})$, $\zeta\Phi(\zeta)>0$ for $\zeta\neq0$, $\Phi(0)=0$, $\Phi'(\zeta)\ge c|\zeta|^{\sigma-1}$, $\sigma\ge1$, for $|\zeta|\le M$.

Our present interest is only to show the applicability of the method. Therefore we assume we are given a weak bounded solution to problem~\eqref{NLNLHE2}. Existence of such a solution is obtained in~\cite{dePabloQuirosRodriguezVazquez2012} in the case $\lj=(-\Delta)^{\alpha/2}$ and $\Phi(\zeta)=|\zeta|^{\sigma-1}\zeta$, $\sigma>(N-\alpha)_+/N$ and quite general initial values. It would be interesting to obtain such a result for general Lévy operators $\lj$.

We prove the nonlinear analogue to Theorem~\ref{th:decay-asymp-sharp} for solution to problem~\eqref{NLNLHE2}. Some restrictions on the exponents make the result not completely satisfactory.

\begin{theorem}
  \label{th:asymp-nonlinear} Assume the hypotheses on the operator $\lj$ of Theorem~{\rm\ref{th:decay-asymp-sharp}} hold and let $u$ be a solution to problem~\eqref{NLNLHE2} with initial value $u_0\in L^1(\mathbb{R}^N)\cap L^\infty(\mathbb{R}^N)$ and $\|u_0\|_\infty\le M$.  Then for every
$\sigma-1<q<p$ (always $q\ge1$), it holds
\begin{equation*}
  \label{decvay-nonl}
  \|u(t)\|_p\le ct^{-\varrho}\|u_0\|_q^\varepsilon
\end{equation*}
for $t$ large, where
\begin{equation}
  \label{exp-decay-nonl}
  \varrho=\varrho(q,p)=\frac{N(p-q)}{p[N(\sigma-1)+\alpha q]},\quad \varepsilon=\varepsilon (q,p)=1-(\sigma-1)\varrho.
\end{equation}
\end{theorem}
The exponents in \eqref{exp-decay-nonl}  are optimal, see \cite{dePabloQuirosRodriguezVazquez2012} for the fractional Laplacian case and $\Phi$ a power.

Following the same idea as before, if we multiply the equation in \eqref{NLNLHE2} by $|u(t)|^{p-2}u(t)$ and integrate, we get the expression
\begin{equation}
  \label{diff-eq-nonl-1}
\frac1p\frac d{dt}\|u(t)\|_p^p=-\mathcal{E}(\Phi(u(t)),|u(t)|^{p-2}u(t)).
\end{equation}
In order to estimate this quantity we first prove a generalized Stroock-Varopoulos inequality.
\begin{proposition}
  \label{pro:SV-J}
  If $u$ is such that $F(u),\,G(u),\,H(u)\in \mathcal{H}_{\lj}(\mathbb{R}^N)$ and $F'G'\ge(H')^2$ then
  \begin{equation}\label{SV-J}
\mathcal{E}(F(u),G(u)) \ge
\mathcal{E}(H(u),H(u)).
\end{equation}
\end{proposition}

\proof
 We only have to look at the original proof of this result with power-type functions performed in \cite{Varopoulos1985} and replace the calculus inequality for powers used there with the following inequality:
for every $a,b\in\RR$ it holds
\begin{equation*}
  \label{calculus-VS}
  (F(b)-F(a))(G(b)-G(a))\ge(H(b)-H(a))^2.
\end{equation*}
  \label{lem:calculus-VS}
In fact, assume without loss of generality that $b>a$, and consider  the case $G(u)=u$. We have
$$
\begin{aligned}
\dyle |H(b)-H(a)|&=\dyle \left|\int_a^bG'(s)\d s\right|\le(b-a)^{1/2}\left(\int_a^b(H')^2(s)\d s\right)^{1/2} \\
&\leq\dyle (b-a)^{1/2}((F(b)-(F(a))^{1/2}.
\end{aligned}
$$
For a general $G$ apply the inequality to $v=G(u)$.
Now the proof of~\eqref{SV-J} follows directly:
$$
\begin{aligned}
\mathcal{E}(F(u),G(u)) &= \frac12\int_{\RR^N}\int_{\RR^N}(F(u(x))-F(u(y)))(G(u(x))-G(u(y)))J(x-y)\d x\d y \\ &\ge
\frac12\int_{\RR^N}\int_{\RR^N}|H(u(x))-H(u(y))|^2J(x-y)\d x\d y=\mathcal{E}(H(u),H(u)).
\end{aligned}$$
\qed

A proof of \eqref{SV-J} when $\lj=(-\Delta)^{\alpha/2}$ is done in~\cite{dePabloQuirosRodriguezVazquez2012} using an extension technique introduced by \cite{CaffarelliSilvestre2007}, based on the fact that $\lj$ is a power of the Laplacian operator. Indeed~\eqref{SV-J} is also valid by the same technique for powers $\lj^{\alpha/2}$ of operators $\lj$ for which it is true.

\noindent\emph{Proof of Theorem~{\rm\ref{th:asymp-nonlinear}}.} Applying the generalized Stroock-Varopoulos inequality~\eqref{SV-J} to expression~\eqref{diff-eq-nonl-1}, and using the behaviour of $\Phi$, we get
$$ \frac{\d}{\d t}\frac1p\int_{\mathbb{R}^N}|u(t)|^p\le -c\mathcal{E}(|u(t)|^{\frac{p+\sigma-1}2},|u(t)|^{\frac{p+\sigma-1}2}).$$
Therefore, Corollary~\ref{cor-inversa-H} with $z=|u|^{\frac{p+\sigma-1}2}$, $s=\frac{2p}{p+\sigma-1}$, $r=\frac{sq}p$ gives, for $\psi(t)=\|u(t)\|_p^p$,
the differential inequality
$$
\psi'(t)+
c\min\left\{\left(\psi(t)\|u_0\|_q^{-p\theta_1}\right)^{\frac1{1-\theta_1}},
\,\left(\psi(t)\|u_0\|_q^{-p\theta_2}\right)^{\frac1{1-\theta_2}}\right\}\le0.
$$
Observe that $\theta_1>\theta_2$. Hence, and exactly as before, only the term with $\theta_1$ has an influence for large times. This implies the estimate~\eqref{decvay-nonl}, where the restriction $1<r<s\le2$ imposes the restriction $\sigma-1<q<p\le p+\sigma-1<2q$.
The proof ends with an iterative process   based on the estimate, for $p_1<p_2<p_3$,
$$
\|u(t)\|_{p_3}\le ct^{-\varrho(p_2,p_3)}\|u(t/2)\|_{p_2}^{\varepsilon(p_2,p_3)}\le ct^{-\varrho(p_2,p_3)}\left(t^{-\varrho(p_1,p_2)}
\|u(0)\|_{p_1}^{\varepsilon(p_1,p_2)}\right)^{\varepsilon(p_2,p_3)}.
$$
We observe that the  exponents \eqref{exp-decay-nonl} satisfy the recurrent relations
  $$
    \varrho(p_1,p_2)\varepsilon(p_2,p_3)+\varrho(p_2,p_3)=\varrho(p_1,p_3), \qquad
    \varepsilon(p_1,p_2)\varepsilon(p_2,p_3)=\varepsilon(p_1,p_3).
  $$
We omit further details.~\qed

%%%%%%%%%%%%%%%%%%%%%%%%%%%%%%%%%%%%%%%%%%%%%%%%%%%%%%%%%%%%%%%%%%%

%%%%%%%%%%%%%%%%%%%%%%%%%%%%%%%%%%%%%%%%%%%%%%%%%%%%%%%%%%%%%%%%%%%

\vskip 1cm

%\newpage
\providecommand{\href}[2]{#2}
\providecommand{\arxiv}[1]{\href{http://arxiv.org/abs/#1}{arXiv:#1}}
\providecommand{\url}[1]{\texttt{#1}}
\providecommand{\urlprefix}{URL }

\end{document}